\documentclass[12pt]{article}

\usepackage{amssymb}
\usepackage{amsmath}
\usepackage[all]{xypic}
\usepackage[frenchb]{babel}
\usepackage{a4wide}

\title{Deux remarques sur le probl\`eme de Lehmer sur les vari\'et\'es ab\'eliennes}

\author{Nicolas Ratazzi}
\date{}

\begin{document}

\maketitle

%\nocite{*}

\renewcommand{\thefootnote}{}
\footnote{\textit{Email address : }ratazzi@math.jussieu.fr}
\renewcommand{\thefootnote}{1}

\vspace{-1cm}

\newcounter{ndefinition}[section]
\newcommand{\defi}{\addtocounter{ndefinition}{1}{\noindent \textbf{D{\'e}finition \thesection.\thendefinition.\ }}}
\newcounter{nrem}
\newcommand{\rem}{\addtocounter{nrem}{1}{\noindent \textbf{Remarque \thesection.\thenrem.\ }}}
\newtheorem{lemme}{Lemme}[section]
\newtheorem{conj}{Conjecture}[section]
\newcounter{nex}[section]
\newcommand{\ex}{\addtocounter{nex}{1}{\noindent \textit{Exemple} \thesection.\thenex\ }}
\newtheorem{theme}{Th{\`e}me} [section]
\newtheorem{prop}{Proposition} [section]
\newtheorem{cor}{Corollaire} [section]
\newtheorem{theo}{Th{\'e}or{\`e}me} [section]
\newcommand{\demo}{\noindent \textit{D{\'e}monstration} : }

\noindent\hrulefill

\noindent \textbf{Abstract :} Let $A/K$ be an abelian variety over a number field $K$. We prove in this article that a good lower bound (in terms of the degree $[K(P):K]$) for the N\'eron-Tate height of the points $P$ of infinite order modulo every strict abelian subvarieties of $A$ implies a good lower bound for the height of all the non-torsion points of $A$. In particular when $A$ is of C.M. type, a theorem of David and Hindry enables us to deduce, up to ``log'' factors, an optimal lower bound for the height of the non-torsion points of $A$. In the C.M. type case, this improves the previous result of Masser \cite{lettre}. Using the same theorem of David and Hindry we prove in the second part an optimal lower bound, up to ``log'' factors, for the product of the N\'eron-Tate height of $n$ End$(A)$-linearly independant non-torsion points of a C.M. type abelian variety.

\vspace{.3cm}



\noindent\hrulefill

\vspace{.3cm}

\noindent On montre ici que la premi\`ere partie de la conjecture de Lehmer ab\'elienne (minoration des points engendrant la vari\'et\'e ab\'elienne en terme de l'indice d'obstruction), formul\'ee dans \cite{davidhindry} entra\^ine la seconde partie de cette conjecture (minoration des points non de torsion en fonction du degr\'e du point et de la dimension du plus petit sous-groupe alg\'ebrique contenant le point). De m\^eme pour le r\'esultat non-conjectural, ce qui permet d'am\'eliorer le pr\'ec\'edent meilleur r\'esultat connu, d\^u \`a Masser \cite{lettre}, pour la minoration des points d'ordre infini sur les vari\'et\'es ab\'eliennes de type C.M. Par ailleurs on montre que la conjecture de Lehmer ab\'elienne entra\^ine la conjecture de Lehmer ab\'elienne multihomog\`ene \textit{a priori} plus forte, telles qu'elles sont \'enonc\'ees dans \cite{davidhindry}. On montre \'egalement que toute avanc\'ee en direction de la conjecture de Lehmer entra\^ine une avanc\'ee similaire en direction de la conjecture multihomog\`ene. En utilisant le r\'esultat principal de \cite{davidhindry} on en d\'eduit, en direction de la conjecture multihomog\`ene, une minoration optimale aux puissances de log pr\`es dans le cas des vari\'et\'es ab\'eliennes de type C.M.

\section{Sur la conjecture de Lehmer sur les vari\'et\'es ab\'e\-lien\-nes}

\noindent Rappelons la conjecture de Lehmer ab\'elienne, formul\'ee dans \cite{davidhindry} conjecture 1.4. On note $\delta_L(P)$ l'indice d'obstruction de $P$.

\vspace{.3cm}

\begin{conj}\label{conj1app2}\textnormal{\textbf{(David-Hindry)}} Soient $A/K$ une vari\'et\'e ab\'elienne de dimension $g$ sur un corps de nombres et $L$ un fibr\'e en droites sym\'etrique ample sur $A$. Il existe une constante strictement positive $c(A/K,L)$ telle que pour tout point $P\in A(\overline{K})$ d'ordre infini modulo toute sous-vari\'et\'e ab\'elienne stricte de $A$, on a
\begin{equation}\label{e1}
\widehat{h}_{L}(P)\geq\frac{c(A/K,L)}{\delta_L(P)}.
\end{equation}
\noindent De plus, en terme du degr\'e $D=[K(P):K]$, on a pour tout point $P\in A(\overline{K})$ qui n'est pas de torsion
\begin{equation}\label{e2}
\widehat{h}_{L}(P)\geq\frac{c(A/K,L)}{D^{\frac{1}{g_0}}},
\end{equation}
\noindent o\`u $g_0$ est la dimension du plus petit sous-groupe alg\'ebrique contenant le point $P$.
\end{conj}

\vspace{.3cm}

\noindent En utilisant le th\'eor\`eme  de David et Hindry \cite{davidhindry} on obtient un r\'esultat, optimal aux puissances de $\log$ pr\`es en direction de l'in\'egalit\'e (\ref{e2}) de la conjecture pr\'ec\'edente. 

\vspace{.3cm}

\begin{theo}\label{thamel}Si $A/K$ est de type C.M., alors il existe une constante strictement positive $c(A/K,L)$ telle que pour tout point $P\in A(\overline{K})$ d'ordre infini, on a 
\[\widehat{h}_{L}(P)\geq\frac{c(A/K,L)}{D^{\frac{1}{g_0}}}\left(\log 2D\right)^{-\kappa(g_0)},\]
\noindent o\`u $D=[K(P):K]$, o\`u $g_0$ est la dimension du plus petit sous-groupe alg\'ebrique de $A$ contenant $P$ et o\`u $\kappa(g_0)=\left(2g_0(g_0+1)!\right)^{g_0+2}$.
\end{theo}
\demo C'est une cons\'equence imm\'ediate du corollaire 2 de \cite{ratazzi} appliqu\'e \`a la vari\'et\'e $V=\overline{\{P\}}$ image sch\'ematique de $P$ dans $A$ sur $K$. On peut faire une preuve directe (ce qui permet d'utiliser le r\'esultat principal de \cite{davidhindry} sans avoir \`a faire intervenir en plus leur remarque utilisant l'indice d'obstruction) : on commence par le cas o\`u $A=\prod_{i=1}^nA_i^{r_i}$, les $A_i$ \'etant des vari\'et\'es ab\'eliennes simples deux \`a deux non-isog\`enes et o\`u $L$ est le fibr\'e en droites ample et sym\'etrique associ\'e au plongement 
\[A=\prod_{i=1}^nA_i^{r_i}\hookrightarrow \prod_{i=1}^n\mathbb{P}_{n_i}^{r_i}\overset{\textnormal{Segre}}{\hookrightarrow} \mathbb{P}^N,\]
\noindent les $A_i$ \'etant plong\'ees dans $\mathbb{P}_{n_i}$ par des fibr\'es $L_i$ tr\`es amples et sym\'etriques. On note $G$ le plus petit sous-groupe alg\'ebrique contenant $V$. On note $G^0$ la composante connexe de l'identit\'e de $G$. C'est une sous-vari\'et\'e ab\'elienne de $A$ et elle est donc isog\`ene \`a $B=\prod_{i=1}^nA_i^{s_i}$ o\`u $0\leq s_i\leq r_i$. On note alors $\pi : A \rightarrow B$ une projection naturelle obtenue par oubli de certaines coordonn\'ees, de sorte que $\pi_{\mid G}$ est une isog\'enie. Montrons que l'on est dans les conditions d'application du th\'eor\`eme principal de \cite{davidhindry} en prenant comme vari\'et\'e ab\'elienne $B$ et comme point $\pi(P)$.

\vspace{.3cm}

Si $\pi(P)$ est d'ordre fini modulo une sous-vari\'et\'e ab\'elienne stricte de $B$, en notant $H$ le plus petit sous-groupe alg\'ebrique contenant $\pi(P)$, on a $\textnormal{dim} H<\textnormal{dim} B$. Ainsi $G_1=G\cap\pi^{-1}(H)$ est un sous-groupe alg\'ebrique strict de $G$ (car $\pi_{\mid G}$ est une isog\'enie), contenant $V$. Ceci est absurde.

\vspace{.3cm}

Si $\pi(P)$ est d'ordre fini, comme $\pi$ est une isog\'enie, le point $P$ est aussi d'ordre fini. Ceci est absurde.

\vspace{.3cm}

\noindent Finalement, $\pi(P)$ est un point d'ordre infini modulo toute sous-vari\'et\'e ab\'elienne de $B$. On peut donc appliquer le th\'eor\`eme principal de \cite{davidhindry}. Par ailleurs, la hauteur et le degr\'e sont d\'efinis relativement aux plongements 
\[A=\prod_{i=1}^nA_i^{r_i}\hookrightarrow \prod_{i=1}^n\mathbb{P}_{n_i}^{r_i}\overset{\textnormal{Segre}}{\hookrightarrow} \mathbb{P}^{N_A}\ \ \textnormal{  et }\ \ B=\prod_{i=1}^nA_i^{s_i}\hookrightarrow \prod_{i=1}^n\mathbb{P}_{n_i}^{s_i}\overset{\textnormal{Segre}}{\hookrightarrow} \mathbb{P}^{N_B}.\]
\noindent De plus l'application $\overline{\pi} :  \prod_{i=1}^n\mathbb{P}_{n_i}^{r_i}\rightarrow  \prod_{i=1}^n\mathbb{P}_{n_i}^{s_i}$ est la projection lin\'eaire d\'efinie par oubli de coordonn\'ees. Dans ce cas, et pour ces plongements, on a 
\[\widehat{h}_{M_B}(\pi(P))\leq \widehat{h}_M(P)\ \ \textnormal{ et }\ \ \deg\pi(P)\leq\deg P.\]
\noindent Ceci nous donne
\begin{align*}
\widehat{h}_M(P) 	& \geq \widehat{h}_{M_B}(\pi(P)),\ \ \ \textnormal{d'o\`u par le th\'eor\`eme de \cite{davidhindry},}\\
			& \geq \frac{c(B, M_B)}{\left(\deg\pi(P)\right)^{\frac{1}{g_0}}}\left(\log 2\deg\pi(P)\right)^{-\kappa(g_0)}\\
			& \geq \frac{c(B, M_B)}{\left(\deg P\right)^{\frac{1}{g_0}}}\left(\log 2\deg P\right)^{-\kappa(g_0)}.\\
			& \geq \frac{c'(A, M)}{\left(\deg P\right)^{\frac{1}{g_0}}}\left(\log 2\deg P\right)^{-\kappa(g_0)},
\end{align*}

\vspace{.3cm}

\noindent o\`u on a pris pour $c'(A,M)$ le minimum des $c(B,M_B)$ quand $s_i$ varie dans $[\![0,r_i]\!]$.

Dans le cas g\'en\'eral, la vari\'et\'e ab\'elienne $A$ est donn\'ee avec une isog\'enie $\rho$ vers la vari\'et\'e ab\'elienne $B=\prod_{i=1}^n A_i^{r_i}$. Soit $P$ d'ordre infini de la vari\'et\'e ab\'elienne de $A$. Le point $Q=\rho(P)$ est un point d'ordre infini de la vari\'et\'e ab\'elienne de $B$. Il r\'esulte facilement de la preuve de la proposition 14. de \cite{phi3} qu'il existe $c'(A,L)$ tel que
\[\widehat{h}_L(P)\geq c'(A,L) \widehat{h}_M(Q).\]
\noindent Ainsi en appliquant le r\'esultat pr\'ec\'edent, on en d\'eduit presque l'in\'egalit\'e voulue : il faut encore remplacer le degr\'e $\deg Q$ par $\deg P$. Or $\deg Q\leq \deg P$. Ceci permet de conclure.\hfill$\Box$

\vspace{.3cm}

\noindent Ce r\'esultat am\'eliore le meilleur r\'esultat pr\'ec\'edemment connu, d\^u \`a Masser qui obtient dans \cite{lettre}, pour tout point $P$ d'ordre infini de $A(\overline{K})$ : 
\[\widehat{h}_{L}(P)\geq \frac{c(A/K,L)}{D^2\log 2D}.\]

\vspace{.3cm}

\noindent En faisant la m\^eme preuve et en appliquant la partie (\ref{e1}) de la conjecture \ref{conj1app2} au lieu du th\'eor\`eme de \cite{davidhindry}, on obtient le

\vspace{.3cm}

\begin{cor}\label{co2} La partie (\ref{e1}) de la conjecture \ref{conj1app2} entra\^ine sa partie (\ref{e2}).
\end{cor}

\section{Sur la conjecture de Lehmer mul\-tiho\-mo\-g\`ene sur les vari\'et\'es ab\'eliennes}

\noindent Soit $A/K$ une vari\'et\'e ab\'elienne de dimension $g$. Quitte \`a augmenter un peu $K$ (cf. par exemple \cite{ratazzi} lemme 1), on peut supposer (et on suppose) que tous les endomorphismes de $A$ sont d\'efinis sur $K$. On note $\widehat{h}_{L}$ la hauteur de N\'eron-Tate sur $A(\overline{K})$ associ\'ee \`a un diviseur ample et sym\'etrique $L$. Pour tout entier $n$ on pose $L_n=L^{\boxtimes n}$ fibr\'e en droites sym\'etrique ample sur $A^n$ et on note $\widehat{h}_{L_n}$ la hauteur de N\'eron-Tate associ\'ee. On commence par un lemme.

\vspace{.3cm}

\begin{lemme}Soit $(P_1,\ldots,P_n)$ un point de $A^n(\overline{K})$. On a 
\[\widehat{h}_{L_n}(P_1,\ldots,P_n)=\sum_{i=1}^n\widehat{h}_{L}(P_i).\]
\end{lemme}
\demo C'est une cons\'equence formelle des propri\'et\'es de fonctorialit\'e des hauteurs de Weil et de la d\'efinition de la hauteur de N\'eron-Tate.\hfill $\Box$

\vspace{.3cm}

\noindent En utilisant ce lemme, on d\'emontre le r\'esultat suivant :

\vspace{.3cm}

\begin{theo}\label{theoapp2}Si $A/K$ est de type C.M., alors, pour tout entier $n\in \mathbb{N}$ il existe une constante $c(A/K,L,n)>0$ telle que pour tout point $(P_1,\ldots,P_n)\in A^n(\overline{K})$ d'ordre infini modulo toute sous-vari\'et\'e ab\'elienne stricte de $A^n$, on a :
\[\prod_{i=1}^n\widehat{h}_{L}(P_i)\geq\frac{c(A/K,L,n)}{D^{\frac{1}{g}}}\left(\log 2D\right)^{-n\kappa(g)},\]
\noindent o\`u $D=[K(P_1,\ldots,P_n):K]$.
\end{theo}
\demo Soient $A_1,\ldots,A_n$ des entiers strictement positifs et $Q_1, \ldots, Q_n$ des points de $A(\overline{K})$ tels que pour tout $i$, $P_i=A_iQ_i$. On a 
\[\widehat{h}_{L_n}(Q_1,\ldots,Q_n)=\sum_{i=1}^n\widehat{h}_{L}(Q_i)=\sum_{i=1}^nA_i^{-2}\widehat{h}_{L}(P_i),\]
\noindent et,
\[\left[K(Q_1,\ldots,Q_n):K\right]^{\frac{1}{ng}}\leq\left(A_1^{2g}\times\cdots\times A_n^{2g}D\right)^{\frac{1}{ng}}.\]
\noindent Le th\'eor\`eme de David-Hindry nous donne alors
\[\sum_{i=1}^nA_i^{-2}\widehat{h}_{L}(P_i)\geq \frac{c(A/K,L,n)}{\left(\prod_{i=1}^nA_i^{\frac{2}{n}}\right)D^{\frac{1}{gn}}}\left(\log \left((\prod_{i=1}^nA_i)D\right)\right)^{-\kappa(g)}.\]
\noindent On pose maintenant, pour tout $1\leq i\leq n$,
\[ x_i=\frac{13\widehat{h}(P_i)}{4\min_j\widehat{h}(P_j)}, \text{ et } A_i=[\sqrt{x_i}].\]
\noindent Pour tout $i$, on a $x_i\geq \frac{13}{4}$ et $x_i\geq A_i^{2}\geq \frac{x_i}{3}$. Ainsi,
\[\sum_{i=1}^n A_i^{-2}\widehat{h}_{L}(P_i)\leq \frac{3\times 4}{13}n\min_j\widehat{h}_{L}(P_j), \]
\noindent et
\[\prod_{i=1}^nA_i^{2}\leq\left(\frac{13}{4\min_j\widehat{h}_{L}(P_j)}\right)^n\prod_{i=1}^n\widehat{h}_{L}(P_i).\]
\noindent Donc,
\begin{align*}
\min_j\widehat{h}_{L}(P_j)	& \geq c_{10}(A/K,L,n)\sum_{i=1}^n A_i^{-2}\widehat{h}_{L}(P_i)\\
				& \geq \frac{c_{11}(A/K,L,n)}{\left(\prod_{i=1}^nA_i^{\frac{2}{n}}\right)D^{\frac{1}{gn}}}\left(\log 2D\prod_{i=1}^nA_i\right)^{-\kappa(g)}\\
				& \geq \frac{4c_{11}(A/K,L,n)\min_j\widehat{h}_{L}(P_j)}{13\prod_{i=1}^n\widehat{h}_{L}(P_i)^{\frac{1}{n}}D^{\frac{1}{gn}}}\left(\log 2D\prod_{i=1}^nA_i\right)^{-\kappa(g)}.
\end{align*}
\noindent Par ailleurs, on a la majoration
\[\log \prod_{i=1}^nA_i\leq n\log\left(\frac{13}{2\min_j\widehat{h}_{L}(P_j)}\right)+2\log \prod_{i=1}^n\widehat{h}_{L}(P_i).\]
\noindent Or on peut toujours supposer que les $\widehat{h}_{L}(P_i)$ sont inf\'erieurs \`a $1$, donc,
\[\log \prod_{i=1}^nA_i\leq n\log\left(\frac{13}{2\min_j\widehat{h}_{L}(P_j)}\right).\]
\noindent Ainsi, 
\[\log 2D\prod_{i=1}^nA_i\leq n\log\left(\frac{13D^{\frac{1}{n}}}{2\min_j\widehat{h}_{L}(P_j)}\right).\]
\noindent On en d\'eduit que 
\[\prod_{i=1}^n\widehat{h}_{L}(P_i)^{\frac{1}{n}}\geq \frac{c_1(A/K,n)}{D^{\frac{1}{ng}}}\left(\log\frac{D^{\frac{1}{n}}}{\min_j\widehat{h}_{L}(P_j)}\right)^{-\kappa(g)}.\]
\noindent Le point $(P_1,\ldots,P_n)$ \'etant d'ordre infini modulo toute sous-vari\'et\'e ab\'elien\-ne, les points $P_i$ sont en particulier d'ordre infini sur $A$. Le r\'esultat inconditionnel de Masser sur la minoration de la hauteur des points sur les vari\'et\'es ab\'eliennes, theorem de \cite{masser}, nous donne donc :
\[\log \frac{D^{\frac{1}{n}}}{\min_j\widehat{h}_{L}(P_j)}\leq c_2(A/K,L,n)\log 2D.\]
\noindent Ainsi, on en d\'eduit
\[\prod_{i=1}^n\widehat{h}_{L}(P_i)\geq \frac{c_3(A/K,L,n)}{D^{\frac{1}{g}}}\left(\log 2D\right)^{-n\kappa(g)}\]
\noindent ce qui conclut.\hfill$\Box$

\vspace{.3cm}

\rem Si au lieu de faire appel au th\'eor\`eme 1.5. de \cite{davidhindry} dans la preuve du th\'eor\`eme \ref{theoapp2} on applique la conjecture \ref{conj1app2}, alors on en d\'eduit le r\'esultat suivant :

\vspace{.3cm}

\begin{theo}\label{theo2}Soient $A/K$ une vari\'et\'e ab\'elienne de dimension $g$ sur le corps de nombres $K$ et $L$ un fibr\'e en droites sym\'etrique ample sur $A$. Si la conjecture \ref{conj1app2} est vraie pour $(A/K,L)$ alors, pour tout entier $n\in \mathbb{N}$ il existe une constante $c(A/K,L,n)>0$ telle que pour tout point $(P_1,\ldots,P_n)\in A^n(\overline{K})$ d'ordre infini modulo toute sous-vari\'et\'e ab\'elienne stricte de $A^n$, on a :
\[\prod_{i=1}^n\widehat{h}_{L}(P_i)\geq\frac{c(A/K,L,n)}{D^{\frac{1}{g}}},\] 
\noindent o\`u $D=[K(P_1,\ldots,P_n):K]$.
\end{theo}

\vspace{.3cm}

\rem En fait dans leur article \cite{davidhindry}, les auteurs formulent \'egale\-ment une conjecture multihomog\`ene du probl\`eme de Lehmer ab\'elien. Plut\^ot que de supposer le point $(P_1,\ldots, P_n)$ d'ordre infini modulo toute sous-vari\'et\'e ab\'elienne stricte de $A^n$, il supposent les points $P_i$ lin\'eairement ind\'ependants dans $A$. Pr\'ecis\'ement ils donnent la conjecture 1.6 suivante :

\vspace{.3cm}

\begin{conj}\label{conj2app2}\textnormal{\textbf{(David-Hindry)}} Soient $A/K$ une vari\'et\'e ab\'elienne de dimension $g$ sur un corps de nombres et $L$ un fibr\'e en droites sym\'etrique ample sur $A$. Pour tout entier $n\in \mathbb{N}$ il existe une constante $c(A/K,L,n)>0$ telle que pour tout $n$-uplet $(P_1,\ldots,P_n)$ de points d'ordre infini dans $A(\overline{K})$, End($A$)-lin\'eairement ind\'ependants, on a :
\[\prod_{i=1}^n\widehat{h}_{L}(P_i)\geq\frac{c(A/K,L,n)}{D^{\frac{1}{g}}},\] 
\noindent o\`u $D=[K(P_1,\ldots,P_n):K]$.
\end{conj}

\vspace{.3cm}

\noindent Dans la formulation de la conjecture \ref{conj2app2} qu'ils donnent, David-Hindry \'ecri\-vent ``li\-n\'e\-aire\-ment ind\'ependants'' sans pr\'eciser s'il s'agit de $\mathbb{Z}$-lin\'eairement ou de End($A$)-lin\'eairement ind\'ependants. Il para\^it pr\'ef\'erable de pr\'eciser. En effet, si on comprend l'assertion ``li\-n\'eaire\-ment ind\'ependants'' comme $\mathbb{Z}$-lin\'eaire\-ment ind\'e\-pen\-dants, alors la conjecture \ref{conj2app2} est fausse comme le montre l'exemple suivant : on prend $E/K$ une courbe elliptique \`a multiplication complexe par un corps quadratique imaginaire contenu dans $K$. On se donne $\alpha\in \textnormal{End}(E)$ un endomorphisme qui n'est pas la multiplication par un entier, on se donne \'egalement un point $P_1$ d'ordre infini dans $E(\overline{K})$ et pour tout $n\geq 1$, on choisit des points $P_n$ tels que $nP_n=P_1$. Enfin on pose $Q_n=\alpha(P_n)$. Puisque $P_1$ est d'ordre infini, les points $P_n$ et $Q_n$ sont $\mathbb{Z}$-lin\'eairement ind\'ependants. De plus on a 
\[\widehat{h}(P_n)\widehat{h}(Q_n)=\frac{\textnormal{N}(\alpha)}{n^4}\widehat{h}(P_1)^2,\]
\noindent et, 
\[D_n:=[K(P_n,Q_n):K]=[K(P_n):K]\leq cn^2.\]
\noindent Donc,
\[\widehat{h}(P_n)\widehat{h}(Q_n)\leq \frac{c'}{D_n^2}.\]
\noindent Ceci montre que l'hypoth\`ese ``$\mathbb{Z}$-lin\'eairement ind\'ependants'' est insuffisante.

\vspace{.3cm}

\noindent Par contre en supposant les points  End($A$)-lin\'eairement ind\'ependants, la situation est bien meilleure. Pr\'ecis\'ement, on a le

\vspace{.3cm}

\begin{theo}\label{theo3}La conjecture \ref{conj1app2} entra\^ine la conjecture \ref{conj2app2}.
\end{theo}
\demo Soit $n>0$ un entier. Au vu du th\'eor\`eme \ref{theo2}, la seule chose \`a prouver, est de montrer que l'hypoth\`ese (i) : ``les points $(P_1,\ldots,P_n)$ sont End($A$)-lin\'eairement ind\'ependants'', entra\^ine l'hypoth\`ese (ii) : ``le point $\mathbf{P}=(P_1,\ldots,P_n)$ est d'ordre infini modulo toute sous-vari\'et\'e ab\'elienne stricte de $A^n$.'' On va plut\^ot montrer que non(ii) implique non(i). Si non(ii) est vraie, alors, il existe un endomorphisme $\varphi$, non-nul, de $A^n$ tel que $\varphi(\mathbf{P})=0$. Or on peut \'ecrire $\varphi(\mathbf{P})=\left(\varphi_1(\mathbf{P}),\ldots,\varphi_n(\mathbf{P})\right)$, o\`u les $\varphi_i$ sont des morphismes de $A^n$ vers $A$ non tous nuls. On suppose par exemple que $\varphi_1$ est non-nul. En notant $\psi_i$ la restriction de $\varphi_1$ \`a la $i$-eme composante de $A^n$, on obtient ainsi $n$ endomorphismes de $A$, $\psi_1,\ldots,\psi_n$, non tous nuls et tels que 
\[\sum_{i=1}^n\psi_i(P_i)=\varphi_1(\mathbf{P})=0.\] 
\noindent Autrement dit, les points $P_1,\ldots,P_n$ sont End($A$)-lin\'eairement d\'ependants.\hfill$\Box$

\vspace{.3cm}

\noindent Enfin la m\^eme preuve permet de constater que le th\'eor\`eme \ref{theo2} entra\^ine un \'enonc\'e analogue en rempla\c{c}ant l'hypoth\`ese ``d'ordre infini modulo toute sous-vari\'et\'e ab\'elienne stricte'' par ``End($A$)-lin\'eairement ind\'ependants''. Ce dernier r\'esultat \`a \'egalement \'et\'e montr\'e par Viada \cite{viada} proposition 4. dans le cas particulier o\`u $A$ est une courbe elliptique.

\vspace{1cm}

\noindent \textbf{Adress :} RATAZZI Nicolas

Universit\'e Paris 6
Institut de Math\'ematiques

Projet Th\'eorie des nombres 

Case 247 

4, place Jussieu 

75252 Paris Cedex 05

FRANCE 

email : ratazzi@math.jussieu.fr

\end{document}